\documentclass[11pt,leqno,final]{amsart}
\usepackage{amssymb,verbatim,enumerate,ifthen}
\usepackage{showkeys}
\usepackage[mathscr]{eucal}
\usepackage[utf8]{inputenc}
\usepackage[T1]{fontenc}

\usepackage{tabu}

\usepackage{dsfont}

\def\N{\mathbb{N}}

\def\R{\mathbb{R}}

\def\D{\mathscr{D}}
\def\W{\mathscr{W}}
\def\Z{\mathscr{Z}}

\def\conv{\mathop{\mbox{\rm conv}}\nolimits}
\def\sgn{\mathop{\text{\rm sgn}}\nolimits}
\def\cl#1{\overline{#1}}
\newtheorem{theorem}{Theorem}
\newtheorem*{theorem*}{Theorem}
\long\def\Thm#1#2{\ifthenelse{\equal{#1}{*}}{\begin{theorem*}#2\end{theorem*}}
             {\begin{theorem}\label{T#1}#2\end{theorem}}}
\newtheorem{Atheorem}{Theorem}

\def\thm#1{Theorem~\ref{T#1}}

\newtheorem{proposition}[theorem]{Proposition}
\newtheorem*{proposition*}{Proposition}
\long\def\Prp#1#2{\ifthenelse{\equal{#1}{*}}{\begin{proposition*}#2\end{proposition*}}
             {\begin{proposition}\label{P#1}#2\end{proposition}}}
\def\prp#1{Proposition~\ref{P#1}}

\newtheorem{corollary}[theorem]{Corollary}
\newtheorem*{corollary*}{Corollary}
\long\def\Cor#1#2{\ifthenelse{\equal{#1}{*}}{\begin{corollary*}#2\end{corollary*}}
             {\begin{corollary}\label{C#1}#2\end{corollary}}}

\newtheorem{lemma}[theorem]{Lemma}
\newtheorem*{lemma*}{Lemma}
\long\def\Lem#1#2{\ifthenelse{\equal{#1}{*}}{\begin{lemma*}#2\end{lemma*}}
             {\begin{lemma}\label{L#1}#2\end{lemma}}}
\def\lem#1{Lemma~\ref{L#1}}

\theoremstyle{definition}
\newtheorem{definition}[theorem]{Definition}
\newtheorem*{definition*}{Definition}
\long\def\Defn#1#2{\ifthenelse{\equal{#1}{*}}{\begin{definition*}\rm #2\end{definition*}}
             {\begin{definition}\label{D#1}\rm #2\end{definition}}}

\newtheorem{remark}[theorem]{Remark}
\newtheorem*{remark*}{Remark}
\long\def\Rem#1#2{\ifthenelse{\equal{#1}{*}}{\begin{remark*}\rm #2\end{remark*}}
             {\begin{remark}\label{R#1}\rm #2\end{remark}}}

\newtheorem{example}{Example}
\newtheorem*{example*}{Example}
\long\def\Exa#1#2{\ifthenelse{\equal{#1}{*}}{\begin{example*}\rm #2\end{example*}}
             {\begin{example}\label{Ex#1}\rm #2\end{example}}}

\def\eq#1{{\rm(\ref{E#1})}}
\def\Eq#1#2{\ifthenelse{\equal{#1}{*}}
  {\begin{equation*}\begin{aligned}[]#2\end{aligned}\end{equation*}}
  {\begin{equation}\begin{aligned}\label{E#1}#2\end{aligned}\end{equation}}}

\begin{document}
%\begin{flushright}
%\emph{Submitted to: Durmogó Dimitriu}
%\end{flushright}
%\vspace{5mm}

\date{\today}

\title{
On a functional equation related to 
two-variable weighted quasi-arithmetic means	
}

\author[T. Kiss]{Tibor Kiss}
\author[Zs. P\'ales]{Zsolt P\'ales}
\address{Institute of Mathematics, University of Debrecen,
H-4032 Debrecen, Egyetem t\'er 1, Hungary}
\email{\{kiss.tibor,pales\}@science.unideb.hu}


\thanks{
The research of the first author has been supported through the new
national excellence program of the ministry of human capacities.	
The research of the second author has been supported by the Hungarian
Scientific Research Fund (OTKA) Grant K111651.}

\begin{abstract}
In this paper, we are going to describe the solutions of the functional equation
\Eq{*}{
\varphi\Big(\frac{x+y}{2}\Big)(f(x)+f(y))
=\varphi(x)f(x)+\varphi(y)f(y)}
concerning the unknown functions $\varphi$ and $f$ defined on an open interval. In our main result only 
the continuity of the function $\varphi$ and a regularity property of the set of zeroes of $f$ are assumed. 
As application, we determine the solutions of the functional equation
\Eq{*}{
  G(g(u)-g(v))=H(h(u)+h(v))+F(u)+F(v)
}
under monotonicity and differentiability conditions on the unknown functions $F,G,H,g,h$.
\end{abstract}

\maketitle

\section{Introduction}

The theory of means have provided a rich source and background for the introduction and investigation of 
functional equations in several variables. The characterization of quasi-arithmetic means was solved 
independently by de Finetti \cite{Def31}, Kolmogorov \cite{Kol30}, and Nagumo \cite{Nag30} for the
case when the number of variables is non-fixed. For the two-variable case, Acz\'el \cite{Acz47b}, 
\cite{Acz48a}, \cite{Acz48b}, \cite{Acz56c}, proved a characterization theorem involving the notion of 
bisymmetry. This result was extended to the $n$-variable case by Maksa--M\"unnich--Mokken
\cite{MunMakMok99}, \cite{MunMakMok00}. A recent characterization theorem of generalized quasi-arithmetic 
means hev been obtained by Matkowski--P\'ales \cite{MatPal15}. Characterization theorems for quasideviation 
means and for Bajraktarevi\'c means were obtained by the author \cite{Pal82a}, \cite{Pal84b}, \cite{Pal87d}.

The equality problem and the so-called invariance equation in various classes of means have been 
investigated in the papers Acz\'el--Kuczma \cite{AczKuc89}, Baj\'ak--P\'ales \cite{BajPal09a}, Berrone 
\cite{Ber00}, Berrone--Lombardi \cite{BerLom01}, Dar\'oczy--Maksa--P\'ales \cite{DarMakPal04}, 
Dar\'oczy--P\'ales \cite{DarPal01a}, \cite{DarPal02c}, Jarczyk \cite{Jar07}, \cite{Jar09b}, \cite{Jar09a}, 
\cite{Jar10a}, \cite{Jar10b}, \cite{Jar13c}, Jarczyk--Matkowski \cite{JarMat06}, \cite{JarMat10}, 
Kahlig--Matkowski \cite{KahMat16}, Leonetti--Matkowski--Tringali \cite{LeoMatTri16}, Losonczi \cite{Los00a}, 
\cite{Los03a}, \cite{Los06b},
Losonczi--P\'ales \cite{LosPal11a}, Mak\'o--P\'ales \cite{MakPal08}, Matkowski \cite{Mat99a}, 
\cite{Mat04c}, \cite{Mat11d}, \cite{Mat14}, Matkowski--Nowicka--Witkowski \cite{MatNowWit17}, P\'ales 
\cite{Pal11}.

In order to solve the equality problem of two-variable functionally weighted quasi-arithmetic means
(called Bajraktarevi\'c means \cite{Baj58}) and quasi-arithmetic means, Z. Daróczy, Gy. Maksa and Zs. 
Páles \cite{DarMakPal04} investigated and solved the functional equation 
\Eq{1}{
\varphi\Big(\frac{x+y}{2}\Big)(f(x)+f(y))
=\varphi(x)f(x)+\varphi(y)f(y),
\qquad(x,y\in I),
}
concerning the unknown functions $\varphi:I\to\R$ and $f:I\to\R$. (Here and throughout this paper, let 
$I\subseteq\R$ stand for a nonempty open interval.) They determined the solutions 
under the natural conditions needed for the definition of weighted quasi-arithmetic means, 
namely, they assumed that $\varphi$ is strictly monotone and continuous, furthermore, that $f$ 
is positive on its domain. The main idea of their approach was to prove that such solutions of \eq{1} are 
infinitely many times differentiable. Thus, using differentiability, and the results of Losonczi \cite{Los99} 
on the equality of two variable Bajraktarevi\'c means, they determined the above equation.  

Although the equation \eq{1} appeared first related to means, its solutions can be interesting 
in general, without any means in the background. Motivated by this, we are going to solve \eq{1} 
assuming only the continuity of $\varphi$ and a regularity property of the zeroes of $f$. After 
eleminating the non-regular solutions, our approach is parallel to what was followed in \cite{DarMakPal04}: 
first, we are going to improve the regularity properties of the unknown functions $\varphi$ and $f$, then, we 
prove that $f$ and $g:=\varphi f$ are solutions of a second-order homogeneous linear differential equation 
with constant coefficients. Solving this differential equation, the solutions of the functional equation 
\eq{1} is finally obtained. Our main result stated in \thm{m1} is proper generalization of \cite[Theorem 
2]{DarMakPal04} eventhough, for a first glance, the two formulations look very different from each other. 
In the construction of the second-order homogeneous linear differential equation with constant coefficients
we also avoided the direct application of the result of Losonczi \cite{Los99}, instead we followed a 
completely independent argument.

As an application of our results, we will solve the functional equation
\Eq{GHF}{
  G(g(u)-g(v))=H(h(u)+h(v))+F(u)+F(v)
}
under monotonicity and differentiability conditions on the unknown functions $F,G,H,g,h$. It remains an open 
problem to reach the same conclusion as in \thm{Fin} without assuming differentiability of the unknown 
functions. A possible way is to follow the regularity improving methods developed in the papers 
Acz\'el--Maksa--P\'ales \cite{AczMakPal99}, \cite{AczMakPal01}, Gil\'anyi--P\'ales \cite{GilPal01a}, 
J\'arai--Maksa--P\'ales \cite{JarMakPal04}, P\'ales \cite{Pal02a}, \cite{Pal03b}.

\section{Auxiliary results}

In the sequel, denote by $\Z_f$ the set of zeros of the function $f$, that is, let
\Eq{*}{\Z_f:=\{x\in I\mid f(x)=0\}=f^{-1}(\{0\}).}
For a subset $S\subseteq I$, we shall also use the notation $S^c:=I\setminus S$. Similarly, $\cl{S}$ will 
denote the closure of the set $S$ relative to $I$. In the next theorem we give a sufficient condition for the 
solutions of equation \eq{1}.

\Thm{sc}{
Let $\varphi,f:I\to\R$ be functions such that $\varphi$ is constant on the set 
$\frac12({\Z_f}^c+I)$. Then the pair $(\varphi,f)$ solves equation \eq{1}.
}

\begin{proof}
Let $x,y\in I$ be arbitrary. Then, obviously, we can distinguish the following four cases. 
\begin{enumerate}[(i)]\itemsep=1mm
\item If $x,y\in \Z_f$, then $f(x)=f(y)=0$ and hence \eq{1} holds trivially.
\item If $x\in\Z_f$ and $y\in{\Z_f}^c$, then $y,\frac{y+x}2\in\frac12({\Z_f}^c+I)$ thus 
$\varphi(y)=\varphi(\frac{x+y}2)$ and $f(x)=0$. These properties imply \eq{1} directly. 
\item The case $y\in\Z_f$ and $x\in{\Z_f}^c$ is analogous to that of (ii).
\item Finally, if $x,y\in{\Z_f}^c$, then $x,y,\frac{y+x}2\in\frac12({\Z_f}^c+I)$ whence 
$\varphi(x)=\varphi(y)=\varphi(\frac{x+y}2)$ follows implying the validity of \eq{1}.
\end{enumerate}
This completes the proof of the statement.
\end{proof}

In order to have a computational formula for the set $\frac12({\Z_f}^c+I)$, we need the 
following lemma.

\Lem{int}{
For any subset $A\subseteq I$, the sum $A+I$ is an open interval, furthermore, if $A\neq\emptyset$, then
\Eq{int}{
A+I=\,]\inf A+\inf I,\sup A+\sup I[\,.
}}

\begin{proof}
If $A\subseteq I$ is empty, then $A+I$ is also empty. Now, assume that $A$ is nonempty and let 
$x,y\in A+I$ be arbitrary with $x<y$. Then there exist $u,v\in A$ such that $x\in u+I$ and $y\in 
v+I$. Then, obviously, $u+v\in(u+I)\cap(v+I)$, furthermore, the sets $u+I$ and $v+I$ are intervals. 
Therefore, if either $x<u+v$ or $u+v<y$, we have either
\Eq{*}{
[x,u+v]\subseteq u+I 
\qquad\mbox{or}\qquad
[u+v,y]\subseteq v+I,
}
respectively. Using these inclusions, if $u+v\leq x$, we obtain that 
\Eq{*}{
[x,y]\subseteq[u+v,y]\subseteq v+I\subseteq A+I.
}
Similarly, if $y\leq u+v$, then
\Eq{*}{
[x,y]\subseteq[x,u+v]\subseteq u+I\subseteq A+I.
}
Finally, if $x<u+v<y$, we get
\Eq{*}{
[x,y]=[x,u+v]\cup[u+v,y]
\subseteq(u+I)\cup(v+I)\subseteq A+I.
}
Hence, the inclusion $[x,y]\subseteq A+I$ holds in each of the above cases, proving that $A+I$ is 
an interval. The openness of $A+I$ is a consequence of the openness of $I$. This directly implies \eq{int}.
\end{proof}

\Lem{?}{
Let $A,B\subseteq I$ be subsets such that $A\subseteq B\subseteq\cl{\conv}(A)$ holds. Then $A+I=B+I$.
}

\begin{proof}
If $A$ is empty, then $B$ is also empty and there is nothing to prove. Therefore, we may restrict ourselves 
to the case when $A$ is not empty. Then, in view of formula \eq{int} of \lem{int}, it is sufficient to show 
that $\inf A=\inf B$ and $\sup A=\sup B$ hold. 

By the condition $A\subseteq B$, it follows that $\inf B\leq\inf A$.
On the other hand, it is easy to see that $\inf\cl{\conv}(A)=\inf A$, therefore the inclusion 
$B\subseteq\cl{\conv}(A)$ implies that $\inf B\geq\inf\cl{\conv}(A)=\inf A$, which completes the proof of 
$\inf A=\inf B$. The proof of $\sup A=\sup B$ is analogous.
\end{proof}

In the next lemma, we establish a consequence of functional equation \eq{1}.

\Lem{0}{
If a pair of functions $(\varphi, f)$ solves \eq{1}, then, for all 
$(x,y)\in\Z_f\times{\Z_f}^c$, we have
\Eq{0}{
\varphi(y)=\varphi\Big(\frac{x+y}2\Big).
}
In addition, if $\varphi$ is continuous, then \eq{0} also holds for all 
$(x,y)\in\cl{\Z}_f\times{\Z_f}^c$.}

\begin{proof}
If $(x,y)\in\Z_f\times{\Z_f}^c$, then $f(x)=0$ and $f(y)\neq0$, hence \eq{1} reduces to \eq{0}. In the case the 
continuity of $\varphi$, the validity of \eq{0} for $(x,y)\in\cl{\Z}_f\times{\Z_f}^c$ follows by a standard 
limiting argument.
\end{proof}

In the following proposition, we give a necessary condition for the solution pairs of equation 
\eq{1}.

\Prp{nc}{
If the pair of functions $(\varphi,f)$ solves \eq{1}, $\varphi$ is continuous on $I$ and $\Z_f$ is 
nonempty, then $\varphi$ is constant on the set $\frac12({\cl{\Z}_f}^c+I)$.}

\begin{proof} The statement is trivial if ${\cl{\Z}_f}^c=\emptyset$, thus we may assume that 
${\cl{\Z}_f}^c$ is nonempty.

The set ${\cl{\Z}_f}^c$ is open, thus it can be written as the union of its components, that 
is, there exists a nonempty set $\Gamma$ and, for all $\gamma\in\Gamma$, there exist extended real 
numbers $a_\gamma<b_\gamma$ such that 
\Eq{*}{{\cl{\Z}_f}^c=\bigcup_{\gamma\in\Gamma}\,]a_\gamma,b_\gamma[\,.}
In the first step we prove that $\varphi$ is constant on any component of ${\cl{\Z}_f}^c$. Let 
$J_\gamma:=\,]a_\gamma,b_\gamma[$ be an arbitrary component of ${\cl{\Z}_f}^c$. Because of that 
$\cl{\Z}_f$ is nonempty, one of the endpoints of $J_\gamma$ must belong to $I$. Without 
loss of generality, we may assume that $a_\gamma\in I$. Then, the maximality of $J_\gamma$ implies 
that $a_\gamma\in\cl{\Z}_f$. Now, let $u\in J_\gamma$ be arbitrarily fixed and define the sequence 
$(u_n)\subseteq J_\gamma$ as $u_n:=2^{1-n}u+(1-2^{1-n})a_\gamma$ whenever $n\in\N$. It is easy to 
see, that $u_1=u$, $u_n\to a_\gamma$ as $n\to\infty$, and that the recursive formula 
$u_{n+1}=\frac{a_\gamma+u_n}{2}$ holds. Now, we are going to show that
\Eq{u}{
\varphi(u)=\varphi(u_n),\qquad(n\in\N).
}
This is trivial for $n=1$. Now, assume that \eq{u} holds for some $n=k\in\N$. To prove \eq{u} for 
$n=k+1$, apply \eq{0} for $y:=u_k\in{\cl{\Z}_f}^c$ and $x:=a_\gamma\in{\cl{\Z}_f}$. Then it follows 
that $\varphi(u_k)=\varphi(\frac{a_\gamma+u_k}{2})=\varphi(u_{k+1})$. By the 
inductive hypothesis, this means that $\varphi(u)=\varphi(u_{k+1})$, which is the required identity.

Upon taking the limit $n\to\infty$ in equation \eq{u} and using the continuity of $\varphi$ at 
$a_\gamma$, it follows that $\varphi(u)=\varphi(a_\gamma)$. Since the point $u\in J_\gamma$ was 
arbitrary, it follows that $\varphi$ is constant on the component $J_\gamma$. 

Now, we show that $\varphi$ is constant on the set $\frac{1}{2}({\cl{\Z}_f}^c+{\cl{\Z}_f}^c)$.
It is easy to check that
\Eq{*}{
\frac{1}{2}({\cl{\Z}_f}^c+{\cl{\Z}_f}^c)=\bigcup_{\gamma,\,\kappa\,\in\,\Gamma} 
\Big]\frac{a_\gamma+a_{\kappa}}{2},\frac{b_\gamma+b_{\kappa}}{2}\Big[\,.
}
If ${\cl{\Z}_f}^c$ has only one component, then we are done. Therefore we may assume that it has at least 
two components, say $J_\gamma=\,]a_\gamma,b_\gamma[\,$ and $J_\kappa=\,]a_\kappa,b_\kappa[\,$ and we can also 
assume that $b_\gamma<a_\kappa$. In this case, because of the maximality of the components $J_\gamma$ and 
$J_\kappa$, we have that $b_\gamma,a_\kappa\in\cl{\Z}_f$. Using what we have already proved,
there exist $\alpha_\gamma,\alpha_\kappa\in\R$, such that $\varphi|_{J_\gamma}=\alpha_\gamma$ and 
$\varphi|_{J_\kappa}=\alpha_\kappa$, moreover, according to the continuity of $\varphi$, it follows that 
$\varphi(b_\gamma)=\alpha_\gamma$ and $\varphi(a_\kappa)=\alpha_\kappa$.  Now, applying the equation \eq{0} for 
$x:=a_\kappa\in\cl{\Z}_f$ and $y\in J_\gamma\subseteq{\cl{\Z}_f}^c$, we get that
$\varphi(\frac{a_\kappa+y}{2})=\varphi(y)=\alpha_\gamma$ for all $y\in J_\gamma$. In other words, $f(u)=\alpha_\gamma$ 
for 
$u\in\frac12(J_\gamma+a_\kappa)=\,\big]\frac{a_\gamma+a_{\kappa}}{2},\frac{b_\gamma+a_{\kappa}}{2}\big[\,$.
Similarly, for 
$u\in\frac12(b_\gamma+J_\kappa)=\,\big]\frac{b_\gamma+a_{\kappa}}{2},\frac{b_\gamma+b_{\kappa}}{2}\big[\,$,
we get that $f(u)=\alpha_\kappa$. By the continuity of $\varphi$ at the point $\frac{b_\gamma+a_{\kappa}}{2}$, 
it follows that $\alpha_\gamma=\alpha_\kappa$ and that $\varphi$ equals the constant $\alpha_\gamma=\alpha_\kappa$ on the 
interval $\,\big]\frac{a_\gamma+a_\kappa}{2},\frac{b_\gamma+b_{\kappa}}{2}\big[\,$. Because 
$\gamma,\kappa\in\Gamma$ was arbitrary, it follows that $\alpha_\gamma$ is the same constant for all 
$\gamma\in\Gamma$, that is, there exists $\alpha\in\R$ such that $\varphi(x)=\alpha$ if 
$x\in\frac{1}{2}({\cl{\Z}_f}^c+{\cl{\Z}_f}^c)$.

To complete the proof, we show that $\varphi$ equals $\alpha$ on $\frac12({\cl{\Z}_f}^c+I)$. Let 
\Eq{*}{
v\in 
\frac12({\cl{\Z}_f}^c+I)
=\big(\tfrac12({\cl{\Z}_f}^c+{\cl{\Z}_f}^c)\big)
\cup\big(\tfrac12({\cl{\Z}_f}^c+\cl{\Z}_f)\big)
}
be arbitrary. We may assume that $v\in\frac12({\cl{\Z}_f}^c+\cl{\Z}_f)$. Then there exists 
$y\in{\cl{\Z}_f}^c\subseteq\big(\frac12({\cl{\Z}_f}^c+{\cl{\Z}_f}^c)\big)\cap{\Z_f}^c$ and 
$x\in\cl{\Z}_f$, such that 
$v=\frac12(x+y)$. By the second assertion of \lem{0}, we have $\varphi(v)=\varphi(y)=\alpha$, which 
finishes the proof.
\end{proof}

\section{Improving of the regularity}

\Prp{ncr}{
If the pair of functions $(\varphi, f)$ solves \eq{1}, $\varphi$ is continuous on $I$ and $\varphi$ 
is not constant on the set $\frac12({\cl{\Z}_f}^c+I)$, then $f$ is nowhere zero and continuous, 
$\varphi$ and $f$ are infinitely many times differentiable, the function $\varphi$ is strictly 
monotone on $I$, and there exists a nonzero real constant $\lambda$ such that 
\Eq{la}{
\varphi'f^2=\lambda
}
holds on $I$.}

\begin{proof}
In view of \prp{nc}, if $\varphi$ is not constant on $\frac{1}{2}({\cl{\Z}_f}^c+I)$ then
$\Z_f$ is empty or, equivalently, $f$ is nowhere zero on $I$.

We claim that there is no nonempty subinterval of $I$ where $\varphi$ would be constant. 
Indirectly, 
assume that this is not the case, that is, there exists a maximal proper subinterval 
$J:=\,]a,b[\,$ of $I$ such that, for some $\alpha\in\R$ we have $\varphi(x)=\alpha$ whenever $x\in 
J$. The interval $J$ cannot be equal to $I$, hence one of its endpoints, say $a\in I$, must belong 
to $I$. The continuity of $\varphi$ implies that $\varphi(a)=\alpha$. Let $u\in I$ and $y\in J$ be 
arbitrarily fixed such that $u<a$ and $\frac12(u+y)\in J$. Then applying \eq{1} for $x\in[u,a]$ and 
$y$, then using that $\varphi(\frac{x+y}{2})=\varphi(y)=\alpha$ and $f(x)\neq 0$, 
we immediately obtain that $\varphi(x)=\alpha$. This contradicts the maximality of $J$. 
Therefore $\varphi$ cannot be constant on any nonempty subinterval of $I$.

Now, we are able to show that $f$ is continuous on $I$. In fact, we are going to show that any point of $I$ 
has a neighborhood where $f$ coincides with a proper continuous function. To do this, let $x_0\in I$ 
be arbitrarily fixed. Then there exists $r>0$ such that $[x_0-2r,x_0+2r]\subseteq I$. Based on 
the previous part of the proof, $\varphi$ cannot be constant on the interval $[x_0-r,x_0+r]$. 
Consequently, there exists $u_0\in[x_0-r,x_0+r]$ such that $\varphi(x_0)\neq\varphi(u_0)$. Let 
$y_0:=2u_0-x_0$. Then $y_0\in [x_0-2r,x_0+2r]$, the function 
$x\mapsto\varphi(\frac{x+y_0}{2})-\varphi(x)$ is continuous on $I$, and, by the choice of $u_0$,
it is different from zero at the point $x:=x_0$. Moreover, there exists $\delta>0$ such that this 
function is different from zero also on the entire interval $]x_0-\delta,x_0+\delta[\,\subseteq I$. 
Using this, equation \eq{1} directly implies that
\Eq{f}{
f(x)
=f(y_0)\frac{\varphi(y_0)-\varphi(\frac{x+y_0}{2})}{\varphi(\frac{x+y_0}{2})-\varphi(x)}
,\qquad(x\in\,]x_0-\delta,x_0+\delta[\,).
}
Thus $f$ equals to a continuous function on the neighborhood $]x_0-\delta,x_0+\delta[\,$ of $x_0$, 
therefore, particularly, $f$ is continuous at $x_0$. Because $x_0$ was arbitrarily chosen, it 
follows that $f$ is continuous on $I$. This implies that $f$ is either positive or negative on $I$.

Thereafter we show that $\varphi$ and $f$ are continuously differentiable on $I$. The argument followed 
here is parallel to that of in the paper \cite{DarMakPal04}. Let $I_0\subseteq I$ be any open subinterval and 
$r>0$ such that the endpoints of the intervals $I_0$, $I_0-r$ and $I_0+r$ belong to $I$. Let further $u\in 
I_0$ and $h\in\,]-r,r[$ be arbitrary. Writing $u+h$ and $u-h$ instead of $x$ and $y$ into the equation 
\eq{1}, respectively, we get that 
\Eq{*}{
\varphi(u)(f(u+h)+f(u-h))=\varphi(u+h)f(u+h)+\varphi(u-h)f(u-h)
} 
holds for all $u\in I_0$ and for all $h\in\,]-r,r[\,$. Integrating both sides of the above equation 
on the interval $]0,r[\,$, then using that $f$ is either positive or negative on $I$, 
a standard calculation 
yields that 
\Eq{ph}{
\varphi(u)=\bigg(\int_{u-r}^{u+r}f(t)\,dt\bigg)^{-1}\int_{u-r}^{u+r}\varphi(t)f(t)\,dt,
\qquad(u\in I_0).
}
By the continuity of $\varphi$ and $f$ on $I_0$, equation \eq{ph} implies that $\varphi|_{I_0}$ is 
continuously differentiable. Hence $\varphi$ is continuously differentiable on $I$. Now, by 
\eq{f}, we easily get that $f$ possesses this property on $I$ too.

In the next step we show that $\varphi$ and $f$ are twice continuously differentiable. After 
differentiating \eq{1} with respect to $x$, we get that
\Eq{d1}{
\frac12\varphi'\Big(\frac{x+y}{2}\Big)(f(x)+f(y))
=(\varphi\cdot f)'(x)-\varphi\Big(\frac{x+y}{2}\Big)f'(x),\qquad(x,y\in I).
}
Keeping the definitions and notations of the previous part, substitute $x:=u+h$ and
$y:=u-h$ into the equation \eq{d1}. Integrating the equation so obtained on the interval
$]0,r[\,$, we get that
\Eq{dph}{
\frac12\varphi'(u)\int_{u-r}^{u+r}f(t)\,dt
=\int_u^{u+r}(\varphi\cdot f)'(t)\,dt-\varphi(u)\int_u^{u+r}f'(t)\,dt,
\qquad(u\in I_0)
}
holds. In view of \eq{dph}, it follows that $\varphi'$ is continuously differentiable on $I_0$, and 
hence also on $I$, therefore $\varphi$ is twice continuously differentiable on $I$. Again, due to 
\eq{f}, we obtain the same conclusion for $f$.

Finally, we prove that $\varphi$ and $f$ are infinitely many times differentiable. To do this,
differentiate \eq{d1} with respect to $y$, and then write $x=y$ into the equation so obtained.  We 
get that
\Eq{*}{
\varphi''(x)f(x)+2\varphi'(x)f'(x)=0,\qquad(x\in I).	
} 
Multiplying this equation by $f(x)$, we can deduce that $(\varphi'\cdot f^2)'=0$. Therefore there 
exists a constant $\lambda\in\R$ such that $\varphi'(x)f^2(x)=\lambda$ for all $x\in I$, where 
$\lambda$ cannot be zero, because $\varphi$ is non-constant. Consequently, \eq{la} is valid.

Applying equations \eq{la} and \eq{f} repeatedly, it can be seen that $\varphi$ and $f$ are
indeed infinitely many times differentiable. Moreover, in view of \eq{la}, we also obtained that
$\varphi$ is strictly monotone on $I$.
\end{proof}

\section{The solutions of the equation \eq{1}}

In order to solve \eq{1}, we rewrite it first into an equivalent form. For two given 
functions $f,g:I\to\R$, define the two-variable function $\D_{f,g}$ by
\Eq{*}{
\D_{f,g}(x,y)
:=\det
\begin{pmatrix}
f\Big(\dfrac{x+y}2\Big) & f(x)+f(y)\\[2mm]
g\Big(\dfrac{x+y}2\Big) & g(x)+g(y)
\end{pmatrix},\qquad(x,y\in I).
}

The following lemma establishes an equivalent form of equation \eq{1} in terms of $\D_{f,g}$.

\Lem{phf}{Let $\varphi,f:I\to\R$ such that $f$ is nowhere zero and define $g:I\to\R$ by  
$g(x):=\varphi(x)f(x)$. Then $(\varphi,f)$ solves \eq{1} if and only if
\Eq{fg}{
\D_{f,g}(x,y)=0,
\qquad(x,y\in I).
}}

\begin{proof}
The equivalence of equations \eq{1} and \eq{fg} can be seen by a short and simple calculation. 
\end{proof}

In order to solve the latter equation for the unknown functions $f$ and $g$, we are going to 
differentiate it twice and four times with respect to the variable $x$. So we obtain differential 
equations for $f$ and $g$. To perform the differentiations, we shall need the 
following extension of the Leibniz Product Rule.

\Lem{L}{
If $B:\R^m\times\R^m\to\R$ is bilinear and $F,G:I\to\R^m$ are $n$ times differentiable  
functions, then $H:I\to\R$ defined by $H(x):=B(F(x),G(x))$ is also $n$ times 
differentiable and
\Eq{H}{
 H^{(n)}(x)=\sum\limits_{k=0}^{n}\binom{n}{k} B\big(F^{(k)}(x),G^{(n-k)}(x)\big), 
\qquad(x\in I).
}}

\begin{proof}In the case when $m=1$ and $B(u,v)=uv$, the above rule is exactly the Leibniz Product  
Rule. However, the proof in the more general case can be carried out (using induction with respect 
to $n$) exactly in the same way as for the case $m=1$.
\end{proof}

\Lem{2}{
Let $f,g:I\to\R$ be $n$ times differentiable functions. Then, for all $x,y\in I$, we have
\Eq{*}{
\partial_1^n\D_{f,g}(x,y)=
\sum\limits_{k=0}^{n-1}\frac{1}{2^{k}}\binom{n}{k}\det
\begin{pmatrix}
f^{(k)}\Big(\dfrac{x+y}2\Big) & f^{(n-k)}(x) \\[2mm]
g^{(k)}\Big(\dfrac{x+y}2\Big) & g^{(n-k)}(x)
\end{pmatrix}+
\frac{1}{2^n}\det
\begin{pmatrix}
f^{(n)}\Big(\dfrac{x+y}2\Big) & f(x)+f(y) \\[2mm]
g^{(n)}\Big(\dfrac{x+y}2\Big) & g(x)+g(y)
\end{pmatrix}.
}}

\begin{proof}  
To prove the assertion of the lemma, let $y\in I$ be fixed, define $B:\R^2\times\R^2\to\R$ by
\Eq{B}{
B\big((u_1,u_2),(v_1,v_2)\big):=\det\begin{pmatrix} u_1 & v_1 \\ u_2 & v_2 \end{pmatrix},}
furthermore define $F,G:I\to\R^2$ by
\Eq{*}{
F(x):=\begin{pmatrix} f\Big(\dfrac{x+y}2\Big) \\[2mm] g\Big(\dfrac{x+y}2\Big) \end{pmatrix}
\qquad\text{and}\qquad
G(x):=\begin{pmatrix} f(x)+f(y) \\[2mm] g(x)+g(y) \end{pmatrix}.
}
Observe that, with the above notations, $H(x)=\D_{f,g}(x,y)$ holds and now the equality \eq{H}  
reduces to the identity to be proved.
\end{proof}

Finally, given a pair $(f,g)$ of sufficiently smooth functions on $I$, define their  
\emph{generalized Wronskian} $\W_{f,g}^{k,\ell}:I\to\R$ for $k,\ell\geq0$ by 
\Eq{wr}{
\W_{f,g}^{k,\ell}(x):=\det
\begin{pmatrix}
f^{(k)}(x) & f^{(\ell)}(x) \\[2mm]
g^{(k)}(x) & g^{(\ell)}(x)
\end{pmatrix}.
}
Obviously, due to the basic properties of the determinant, $\W_{f,g}^{k,\ell}=-\W_{f,g}^{\ell,k}$,  
therefore, the function $\W_{f,g}^{k,\ell}$ is identically zero on $I$ if $k=\ell$.

\Thm{dfe}{
Let $f,g:I\to\R$ be 4 times differentiable functions such that $\W_{f,g}^{0,1}$ is not identically  
zero on $I$. Then $(f,g)$ solves \eq{fg} if and only if there exist constants 
$a,b,c,d,\gamma\in\R$ with $ad\neq bc$ such that 
\Eq{abcd}{
\begin{pmatrix}
f \\[1mm] g
\end{pmatrix}=
\begin{pmatrix}
a & b \\ c & d
\end{pmatrix}
\begin{pmatrix}
\psi_1 \\[1mm] \psi_2
\end{pmatrix},
}
where, for all $x\in I$, 
\begin{enumerate}\itemsep=1mm
\item $\psi_1(x)=\sin\big(\sqrt{-\gamma}x\big)$ and $\psi_2(x)=\cos\big(\sqrt{-\gamma}x\big)$
if $\gamma<0$,
\item $\psi_1(x)=x$ and $\psi_2(x)=1$ if $\gamma=0$, and
\item $\psi_1(x)=\sinh\big(\sqrt{\gamma}x\big)$ and $\psi_2(x)=\cosh\big(\sqrt{\gamma}x\big)$ if  
$\gamma>0$.
\end{enumerate}
}

\begin{proof}
Assume that the pair $(f,g)$ solves \eq{fg}. First we are going to show that there exist constants  
$\alpha,\beta\in\R$ such that
\Eq{c}{
\W_{f,g}^{0,1}(x)=\alpha \qquad\text{and}\qquad
\W_{f,g}^{1,2}(x)=\beta,\qquad(x\in I).
}
Let $y\in I$ be fixed. Differentiating \eq{fg} with respect to the first variable twice and 
four times, then, applying \lem{2} for $n=2$ and $n=4$, we get that
\Eq{*}{
\det
\begin{pmatrix}
f\Big(\dfrac{x+y}2\Big) & f''(x) \\[2mm]
g\Big(\dfrac{x+y}2\Big) & g''(x)
\end{pmatrix}
+\det
\begin{pmatrix}
f'\Big(\dfrac{x+y}2\Big) & f'(x) \\[2mm]
g'\Big(\dfrac{x+y}2\Big) & g'(x)
\end{pmatrix}
+\frac{1}{4}\det
\begin{pmatrix}
f''\Big(\dfrac{x+y}2\Big) & f(x)+f(y) \\[2mm]
g''\Big(\dfrac{x+y}2\Big) & g(x)+g(y)
\end{pmatrix}=0
}
and
\Eq{*}{
\det
\begin{pmatrix}
f\Big(\dfrac{x+y}2\Big) & f^{(4)}(x) \\[2mm]
g\Big(\dfrac{x+y}2\Big) & g^{(4)}(x)
\end{pmatrix}
+&2\det
\begin{pmatrix}
f'\Big(\dfrac{x+y}2\Big) & f'''(x) \\[2mm]
g'\Big(\dfrac{x+y}2\Big) & g'''(x)
\end{pmatrix}
+\frac32\det
\begin{pmatrix}
f''\Big(\dfrac{x+y}2\Big) & f''(x) \\[2mm]
g''\Big(\dfrac{x+y}2\Big) & g''(x)
\end{pmatrix}\\
+&\frac12\det
\begin{pmatrix}
f'''\Big(\dfrac{x+y}2\Big) & f'(x) \\[2mm]
g'''\Big(\dfrac{x+y}2\Big) & g'(x)
\end{pmatrix}
+\frac1{16}\det
\begin{pmatrix}
f^{(4)}\Big(\dfrac{x+y}2\Big) & f(x)+f(y) \\[2mm]
g^{(4)}\Big(\dfrac{x+y}2\Big) & g(x)+g(y)
\end{pmatrix}=0
}
holds for all $x\in I$. Substituting $x=y\in I$ into the previous equations, they reduce to
\Eq{ew1}{
\W_{f,g}^{0,2}(x)=0
\qquad\text{and}\qquad
7\W_{f,g}^{0,4}+12\W_{f,g}^{1,3}(x)=0,
}
respectively.

For $k,\ell\geq 0$, define the functions $F_k,G_\ell:I\to\R^2$ as
\Eq{*}{
F_k(x):=
\begin{pmatrix}
f^{(k)}(x) \\ g^{(k)}(x)
\end{pmatrix}
\qquad\text{and}\qquad
G_\ell(x):=
\begin{pmatrix}
f^{(\ell)}(x) \\ g^{(\ell)}(x)
\end{pmatrix},
}
and, in the rest of the proof, let $B:\R^2\times\R^2\to\R$ be defined as in \eq{B}. 

Applying \lem{L} for $n=2$ and for the functions $B$, $F:=F_0$ and $G:=G_2$, we get that the  
identity $\big(\W_{f,g}^{0,2}\big)''=2\W_{f,g}^{1,3}+\W_{f,g}^{0,4}$ holds on $I$, thus the system 
of equations \eq{ew1} is equivalent to the following one:
\Eq{ew2}{
\W_{f,g}^{0,2}(x)=0
\qquad\text{and}\qquad
\W_{f,g}^{1,3}(x)=0.
}
By obvious application of \lem{L} (for $n=1$, $B$, $F:=F_0$, $G:=G_1$ and $n=1$, $B$, $F:=F_1$,  
$G:=G_2$), we obtain that $\big(\W_{f,g}^{0,1}\big)'=\W_{f,g}^{0,2}$ and 
$\big(\W_{f,g}^{1,2}\big)'=\W_{f,g}^{1,3}$ hold on $I$. In view of \eq{ew2}, this means that there 
exist constants $\alpha,\beta\in\R$ such that \eq{c} holds for all $x\in I$. By our assumption 
$\alpha$ is different from zero. Now, consider the equation
\Eq{ps0}{
\det
\begin{pmatrix}
\psi&\psi'&\psi''\\[1mm]
f&f'&f''\\[1mm]
g&g'&g''
\end{pmatrix}
=0
}
on $I$ with the unknown function $\psi:I\to\R$. Expanding the determinant with respect to its first 
row, in view of \eq{c}, we get that \eq{ps0} is a homogeneous second-order linear differential 
equation and it is equivalent to
\Eq{ps1}{
\psi''=\gamma\psi,
}
where $\gamma$ denotes the constant $-\beta/\alpha$.
Then, based on the theory of linear differential equations with constant coefficients, we have that 
the functions $\psi_1$ and $\psi_2$ form a solution basis of \eq{ps1} in each of the possibilities 
(1), (2) or (3). On the other hand, the functions $f$ and $g$ trivially solve \eq{ps0}, and hence 
also \eq{ps1}. Thus, they must be linear combinations of $\psi_1$ and $\psi_2$, that is, \eq{abcd} 
holds for some constants $a,b,c,d\in\R$ with $ad\neq bc$.

To prove the reversed statement, assume now that there exist $a,b,c,d,\gamma\in\R$ with $ad\neq bc$ 
such that \eq{abcd} holds. By the product rule for determinants, we have that
\Eq{*}{
\D_{f,g}(x,y)=\det\begin{pmatrix} a & b \\ c & d \end{pmatrix} \D_{\psi_1,\psi_2}(x,y), 
\qquad(x,y\in I),
}
therefore, to check the validity of \eq{fg}, it is sufficent to prove that 
\Eq{12}{
\D_{\psi_1,\psi_2}(x,y)=0,
\qquad(x,y\in I).
}
This equality is obvious if $\gamma=0$. To see this in the remaining cases, assume first that  
$\gamma>0$. Then, by the well-known identities for hyperbolic functions, we get that
\Eq{*}{
\psi_1(x)+\psi_1(y)
&=\sinh\big(\sqrt{\gamma}x\big)+\sinh\big(\sqrt{\gamma}y\big) \\
&=2\sinh\Big(\sqrt{\gamma}\frac{x+y}{2}\Big)\cosh\Big(\sqrt{\gamma}\frac{x-y}{2}\Big)
=2\psi_1\Big(\frac{x+y}{2}\Big)\cosh\Big(\sqrt{\gamma}\frac{x-y}{2}\Big)
}
and
\Eq{*}{
\psi_2(x)+\psi_2(y)
&=\cosh\big(\sqrt{\gamma}x\big)+\cosh\big(\sqrt{\gamma}y\big) \\
&=2\cosh\Big(\sqrt{\gamma}\frac{x+y}{2}\Big)\cosh\Big(\sqrt{\gamma}\frac{x-y}{2}\Big)
=2\psi_2\Big(\frac{x+y}{2}\Big)\cosh\Big(\sqrt{\gamma}\frac{x-y}{2}\Big)
}
hold for all $x,y\in I$. Therefore, using these formulae, we can see that the columns of the 
determinant 
$\D_{\psi_1,\psi_2}(x,y)$ are linearly dependent, which implies \eq{12}.

The proof of \eq{12} in the case $\gamma<0$ is based on similar identities for trigonometric  
functions, and therefore it is completely analogous.
\end{proof}

\Thm{m1}{
Let $\varphi,f:I\to\R$ such that $\varphi$ is continuous and ${\Z_f}^c\subseteq\cl{\conv}({\cl{\Z}_f}^c)$. 
Then $(\varphi,f)$ solves functional equation \eq{1} if and only if either
\begin{enumerate}[(a)]
\item there exists an interval $J\subseteq I$ such that $f(x)=0$ for all $x\in I\setminus J$  
and $\varphi$ is constant on $\frac12(J+I)$,
\end{enumerate}
or
\begin{enumerate}[(b)]
\item there exist constants $a,b,c,d,\gamma\in\R$ with $ad\neq bc$ such that 
\begin{enumerate}[(1)]
\item $f(x)=a\sin(\sqrt{-\gamma}x)+b\cos(\sqrt{-\gamma}x)\neq0$ and 
$\varphi(x)=\dfrac{c\sin(\sqrt{-\gamma}x)+d\cos(\sqrt{-\gamma}x)}{a\sin(\sqrt{-\gamma}x)+b\cos(\sqrt{-\gamma}
x)}$ if $\gamma<0$,
\item $f(x)=a+bx\neq0$ and $\varphi(x)=\dfrac{c+dx}{a+bx}$ if $\gamma=0$,
\item $f(x)=a\sinh(\sqrt{\gamma}x)+b\cosh(\sqrt{\gamma}x)\neq0$ and 
$\varphi(x)
=\dfrac{c\sinh(\sqrt{\gamma}x)+d\cosh(\sqrt{\gamma}x)}{a\sinh(\sqrt{\gamma}x)+b\cosh(\sqrt{\gamma}x)}$ if 
$\gamma>0$
\end{enumerate}
\end{enumerate}
for all $x\in I$.
}

\begin{proof} 

We can distinguish the following two main cases:
\begin{enumerate}[(i)]\itemsep=1mm
 \item either $\Z_f\neq\emptyset$ or $\Z_f=\emptyset$ and $\varphi$ is constant on $I$, or
 \item $\Z_f=\emptyset$ and $\varphi$ is non-constant on $I$.
\end{enumerate}
Consider first the case (i) and assume that $(\varphi,f)$ solves equation \eq{1}. If ${\Z_f}^c$ is empty, 
then ${\cl{\Z}_f}^c$ is also empty and we trivially have $\frac12({\cl{\Z}_f}^c+I)=\frac12({\Z_f}^c+I)$. If 
${\Z_f}^c$ is nonempty then the condition ${\Z_f}^c\subseteq\cl{\conv}({\cl{\Z}_f}^c)$ implies that the set 
${\cl{\Z}_f}^c$ is also nonempty. Applying \lem{?} for $A:={\cl{\Z}_f}^c$ and $B:={\Z_f}^c$, we obtain again 
that $\frac12({\cl{\Z}_f}^c+I)=\frac12({\Z_f}^c+I)$ hold.

Using this, we show that $\varphi$ is constant on the set $\frac12({\Z_f}^c+I)$. If $\Z_f$ is nonempty, then, 
due to \prp{nc}, the function $\varphi$ is constant on $\frac12({\cl{\Z}_f}^c+I)=\frac12({\Z_f}^c+I)$. 
If $\Z_f$ is empty, then $\frac12({\Z_f}^c+I)=I$, thus, by the second condition of case (i), we also have 
that $\varphi$ is constant on $\frac12({\Z_f}^c+I)$.

Let $J$ denote the convex hull of the set ${\Z_f}^c$. Then $I\setminus J\subseteq I\setminus{\Z_f}^c=\Z_f$,  
which implies 
that $f(x)=0$ whenever $x\in I\setminus J$. Finally, we show that $\varphi$ is constant on $\frac12(J+I)$. If 
${\Z_f}^c$ is empty then $J$ is also empty, hence $\varphi$ is trivially constant on $\frac12(J+I)=\emptyset$. 
In the other case when ${\Z_f}^c$ is nonempty, we have that $\{\inf J,\sup J\}=\{\inf\Z_f^c,\sup\Z_f^c\}$, 
therefore, formula \eq{int} of \lem{int} implies that $\frac12({\Z_f}^c+I)=\frac12(J+I)$. Thus $\varphi$ is 
constant on $\frac12(J+I)$, and consequently, (a) holds.

Conversely, assume that the alternative (a) is valid, that is, there exists an interval $J\subseteq I$ 
such that $I\setminus J\subseteq\Z_f$ and $\varphi$ is constant on $\frac12(J+I)$. Then ${\Z_f}^c\subseteq J$, 
hence $\varphi$ is constant on $\frac12({\Z_f}^c+I)$, which, in view of \thm{sc}, implies that $\varphi$ and 
$f$ solve \eq{1}.

Now, consider case (ii), namely suppose that $\Z_f$ is empty and $\varphi$ is non-constant on $I$. Then, by 
\prp{ncr}, the functions $\varphi,f$ are infinitely many times differentiable such that \eq{la} 
holds with a nonzero constant $\lambda$. Then, in view of \lem{phf}, the functions $f$ and $g:=\varphi f$ are 
infinitely many times differentiable solutions of \eq{fg}. Furthermore, we have that 
$\W_{f,g}^{0,1}=\lambda\neq0$. Applying \thm{dfe}, we obtain that $\varphi=\frac{g}{f}$ and $f$ must be one 
of the forms (1), (2) or (3) represented in the alternative (b).

Assuming (b), one can easily see that, in each cases (1), (2) and (3), the functions $f$ and $g:=\varphi f$
are solutions of \eq{fg}. Consequently, $\varphi$ and $f$ are solutions of \eq{1}.
\end{proof}

\Rem{m1}{The condition ${\Z_f}^c\subseteq\cl{\conv}({\cl{\Z}_f}^c)$ is easily fulfilled if either $f$ is 
nowhere zero (then $\Z_f$ is empty) or $f$ is continuous (then $\Z_f$ is closed). Therefore, \thm{m1} is a 
proper generalization of the result in \cite{DarMakPal04}, where a similar conclusion was reached assuming 
that $\varphi$ was strictly monotone and continuous and $f$ was positive. It seems to be an interesting 
question if the conclusion of \thm{m1} could be obtained without assuming the condition 
${\Z_f}^c\subseteq\cl{\conv}({\cl{\Z}_f}^c)$.}

\section{Application}

In this section we are going to solve the functional equation
\Eq{GHF+}{
G(g(u)-g(v))=H(h(u)+h(v))+F(u)+F(v), \qquad(u,v\in J),
}
where $J\subseteq\R$ is a nonempty open interval, furthermore 
\Eq{*}{
g,h,F:J\to\R,\quad
G:g(J)-g(J)\to\R,\qquad\mbox{and}\qquad
H:h(J)+h(J)\to\R
}
are considered as unknown functions. 

\Lem{GHF+}{With the above notations, assume that $h:J\to\R$ is a continuous strictly monotone function and 
define
\Eq{not1}{
   I:=2h(J),\qquad \ell(x):=g\circ h^{-1}\Big(\frac{x}{2}\Big), \qquad(x\in I), 
   \qquad\mbox{and}\qquad G_0:=G-G(0).
}
If $(g,h,F,G,H)$ is a solution of \eq{GHF+}, then $\ell,H:I\to\R$ and $G_0:\ell(I)-\ell(I)\to\R$ 
solve the functional equation
\Eq{G0H}{
  G_0(\ell(x)-\ell(y))=H\Big(\frac{x+y}{2}\Big)-\frac{H(x)+H(y)}2, \qquad(x,y\in I)
}
and 
\Eq{F}{
F(u)=\frac12\big(G(0)-H(2h(u))\big), \qquad(u\in J).
}
Conversely, if $\ell,H:I\to\R$ and $G_0:\ell(I)-\ell(I)\to\R$ solve the functional equation \eq{G0H}, 
$G(0)\in\R$ is an arbitrary constant, $J:=h^{-1}(\frac12 I)$, $F:J\to\R$ is defined by \eq{F}, furthermore 
$g:J\to\R$ and $G:\ell(I)-\ell(I)\to\R$ are given by 
\Eq{*}{
   g(u):=\ell(2h(u)) \qquad\mbox{and}\qquad G:=G(0)+G_0,
}
then $(g,h,F,G,H)$ is a solution of \eq{GHF+}.}

\begin{proof}
Substituting $v=u$ into \eq{GHF+}, it immediately follows that $F$ is of the form \eq{F}.
Therefore, \eq{GHF+} can be rewritten as
\Eq{GH}{
  G(g(u)-g(v))-G(0)=H(h(u)+h(v))-\frac{H(2h(u))+H(2h(v))}2 \qquad(u,v\in J).
}
Since $h:J\to\R$ is a continuous strictly monotone function, thus its inverse is also a continuous 
strictly monotone function which is defined on the open interval $h(J)$. With the notation \eq{not1}, after 
the substitution $u:=h^{-1}\big(\frac{x}{2}\big)$, $v:=h^{-1}\big(\frac{y}{2}\big)$, equation \eq{GH} reduces 
to \eq{G0H}.

The proof of the reversed implication is a simple computation, therefore, it is omitted.
\end{proof}

\Lem{GH}{Assume that $\ell,H:I\to\R$ and $G_0:\ell(I)-\ell(I)\to\R$ are differentiable solutions of 
functional equation \eq{G0H} such that $\ell'$ does not vanish on $I$. Then the pair of functions 
$(\varphi,f)$ given by
\Eq{not}{
  \varphi:=H' \qquad\mbox{and}\qquad f:=\frac{1}{\ell'}
}
solves equation \eq{1}.}

\begin{proof} Differentiating equation \eq{G0H} with respect to the variables $x$ and $y$, we get
\Eq{*}{
  G_0'(\ell(x)-\ell(y))\cdot\ell'(x)&=\frac12H'\Big(\frac{x+y}{2}\Big)-\frac{H'(x)}2, \\
 -G_0'(\ell(x)-\ell(y))\cdot\ell'(y)&=\frac12H'\Big(\frac{x+y}{2}\Big)-\frac{H'(y)}2
}
for all $x,y\in I$. Multiplying the first equation by $\frac{2}{\ell'(x)}$, the second one by 
$\frac{2}{\ell'(y)}$ and adding up the equations so obtained side by side, we obtain that
\Eq{*}{
0=H'\Big(\frac{x+y}{2}\Big)\Big(\frac{1}{\ell'(x)}+\frac{1}{\ell'(y)}\Big)
-\frac{H'(x)}{\ell'(x)}-\frac{H'(y)}{\ell'(y)}, \qquad(x,y\in I).
}
Therefore, with the notations \eq{not}, we can see that \eq{1} is satisfied.
\end{proof}

Finally we give the complete solution of \eq{G0H}.

\Thm{Fin}{Let $\ell,H:I\to\R$ and $G_0:\ell(I)-\ell(I)\to\R$ are differentiable functions 
such that $\ell'$ does not vanish on $I$ and $H'$ is continuous. Then the triple $(G_0,\ell,H)$ solves 
functional equation \eq{G0H} if and only if there exist constants $A,B,C,D,E,\alpha,\beta$ with 
$CD\alpha\neq0$ such that, for all $x\in I$ and $u\in\ell(I)-\ell(I)$ one of the following possibilities 
holds:
\begin{center}
\begin{tabular}{l|l|l|l}
                & $G_0(u)=$         & $\ell(x)=$                                  & $H(x)=$ \\[.3mm]\hline
\vspace{-3.5mm} &&&\\
(i)             & $0$               & arbitrary                                   & $Ax+B$ \\[.5mm]
(ii)            & $C(Du)^2$         & $\frac1{2D}x+E$                             & $-Cx^2+Ax+B$  \\[.5mm]
(iii)           & $C(Du)^2$         & $\frac1{2D}e^{\alpha x}+E$                  & $-\frac{C}2e^{2\alpha x}+Ax+B$  \\[.5mm]
(iv)            & $C\ln(\cosh(Du))$ & $\frac1{2D}\ln|\alpha x+\beta|+E$           & $C\ln|\alpha x+\beta|+Ax+B$  \\[.5mm]
(v)             & $C\ln(\cosh(Du))$ & $\frac1{2D}\ln|\tan(\alpha x+\beta)|+E$     & $C\ln|\sin(2\alpha x+2\beta)|+Ax+B$ \\[.5mm]
(vi)            & $C\ln(\cosh(Du))$ & $\frac1{2D}\ln|\tanh(\alpha x+\beta)|+E$    & $C\ln|\sinh(2\alpha x+2\beta)|+Ax+B$ \\[.5mm]
(vii)           & $C\ln(\cos(Du))$  & $\frac1{D}\arctan(\tanh(\alpha x+\beta))+E$ & $C\ln(\cosh(2\alpha x+2\beta))+Ax+B$
\end{tabular}
\end{center}
In addition, in cases (iv), (v), and (vi) we have 
\Eq{dom}{
  0\notin \alpha I+\beta, \qquad 
  \mathbb{Z}\cap\frac{2}{\pi}\big(\alpha I+\beta\big)=\emptyset, 
  \qquad\mbox{and}\qquad 0\notin \alpha I+\beta,
}
respectively.}

\begin{proof} In the proof we are going to combine the result of \thm{m1} and \lem{GH}.

Suppose that $(G_0,\ell,H)$ solves \eq{G0H}. Under the assumptions of the theorem, $\varphi=H'$ and 
$f=\frac1{\ell'}$ are solutions of the functional equation \eq{1} and $f$ is nowhere zero. 
Therefore, by \thm{m1}, either $\varphi$ is constant on $I$ or there exist constants 
$a,b,c,d,\gamma$ with $ad\neq bc$ such that $f$ and $\varphi$ are of the forms (1), (2), or (3) in 
alternative (b) of \thm{m1}.

If $\varphi$ is constant, then $H$ is affine, that is, there exist constants $A,B\in\R$ such that 
$H(x)=Ax+B$ for $x\in I$. Then the right hand side of \eq{G0H} is identically zero, therefore, $G_0$ must equal to zero on 
$\ell(I)-\ell(I)$ and $\ell$ can be arbitrary. That is, case \emph{(i)} holds.

From now on assume that $\varphi=H'$ is non-constant on $I$. We distinguish six main cases 
according to the sign of the parameter $\gamma$ and the parameters $a,b$. In each cases we are 
going to solve the differential equations in \eq{not} for the unknown functions $H$ and $\ell$, 
and then we determine $G_0$ using the functional equation \eq{G0H}.

%----------Trigonometrikus
\emph{Case 1. Assume that $\gamma<0$.} By the condition $ad\neq bc$, we have that $a^2+b^2>0$. 
Denote $\alpha:=\frac{\sqrt{-\gamma}}2\neq0$ and define $\beta\in[0,\pi[$ as the unique solution of the system of
equations
\Eq{*}{
\cos(2\beta)=\frac{a}{\sqrt{a^2+b^2}},\qquad 
\sin(2\beta)=\frac{b}{\sqrt{a^2+b^2}}.
}
Then we have that 
\Eq{*}{
 0\neq \frac{1}{\ell'(x)}
 =a\sin(\sqrt{-\gamma}x)+b\cos(\sqrt{-\gamma}x)
 =\sqrt{a^2+b^2}\sin(2\alpha x+2\beta),\qquad(x\in I).
}
Therefore, $2\alpha x+2\beta\notin\pi\mathbb{Z}$ for all $x\in I$, which yields that the second condition in 
\eq{dom} is satisfied.

Solving the differential equations in \eq{not} for the unknown functions $\ell$ and $H$, with the notations
$A:=\frac{ac+bd}{a^2+b^2}$, $C:=\frac{ad-bc}{\sqrt{-\gamma}(a^2+b^2)}\neq 0$, and 
$D:=\frac{\sqrt{-\gamma(a^2+b^2)}}{2}\neq 0$, we get that there exist constants $E$ and $B$ such 
that, for all $x\in I$, we have
\Eq{*}{
\ell(x)=\frac1{2D}\ln|\tan(\alpha x+\beta)|+E\qquad\text{and}\qquad
H(x)=C\ln|\sin(2\alpha x+2\beta)|+Ax+B.
}
Using these representations and substituting $v:=\alpha x+\beta$, $w:=\alpha y+\beta$, the functional equation 
\eq{G0H} reduces to
\Eq{*}{
  G_0\Big(\frac{1}{2D}\ln\frac{\tan v}{\tan w}\Big) 
  &=C\ln|\sin(v+w)|-\frac{C\ln|\sin(2v)|+C\ln|\sin(2w)|}{2} \\
  &=C\ln\frac{|\sin(v+w)|}{\sqrt{\sin(2v)\sin(2w)}}
  =C\ln\frac{|\sin v\cos w+\sin w\cos v|}{\sqrt{4\sin v\cos v\sin w\cos w}} 
  =C\ln\frac{\sqrt{\frac{\tan v}{\tan w}}+\sqrt{\frac{\tan w}{\tan v}}}{2}.
}
(Note that, due to the condition $\mathbb{Z}\cap\frac{2}{\pi}\big(\alpha I+\beta\big)=\emptyset$, in the 
above computation we have that $\tan v$ and $\tan w$ as well as $\sin(2v)$ and $\sin(2w)$ have the same sign.)
Putting $u:=\frac{1}{2D}\ln\frac{\tan v}{\tan w}$, it follows that $\sqrt{\frac{\tan v}{\tan 
w}}=e^{Du}$, 
hence the above equation yields that
\Eq{*}{
G_0(u)=C\ln(\cosh(Du)).
}
That is, we obtain the solutions listed in \emph{(v)}.

%----------Affin(a)
\emph{Case 2. Assume that $\gamma=0$ and $b=0$.} Then, in view of the condition $ad\neq bc$, the 
parameters $a$ and $d$ are different from zero. Solving \eq{not}, with the notations 
$A:=\frac{c}{a}$, $C:=-\frac{d}{2a}\neq0$ and $D:=\frac{a}{2}\neq0$, we obtain that there exist 
constants $E$ and $B$ such that, for all $x\in I$, we have
\Eq{*}{
\ell(x)=\frac1{2D} x+E\qquad\text{and}\qquad
H(x)=-Cx^2+Ax+B.
}
Then the equation \eq{G0H} reduces to the form
\Eq{*}{
  G_0\Big(\frac{x-y}{2D}\Big)=-C\Big(\frac{x+y}{2}\Big)^2+\frac{Cx^2+Cy^2}{2}
  =C\Big(\frac{x-y}{2}\Big)^2.
}
Now, replacing $\frac{x-y}{2D}$ by $u$, it follows that $G_0(u)=C(Du)^2$. That is, the solutions 
obtained are exactly those listed in case \emph{(ii)}. 

%----------Affin(b)
\emph{Case 3. Assume that $\gamma=0$ and $b\neq0$.} Solving the equations in \eq{not} and using 
the notations $A:=\frac{d}{b}$, $C:=\frac{bc-ad}{b^2}\neq0$, $D:=\frac{b}{2}\neq0$, $\alpha:=b\neq0$ 
and $\beta:=a$, we can see that $0\notin \alpha I+\beta$ (that is, the first condition in \eq{dom} holds) and 
there exist constants $E$ and $B$ such that, for all $x\in I$,
\Eq{*}{
\ell(x)=\frac1{2D}\ln|\alpha x+\beta|+E\qquad\text{and}\qquad
H(x)=C\ln|\alpha x+\beta|+Ax+B
}
holds. Consequently, introducing the notations $v:=|\alpha x+\beta|$ and $w:=|\alpha y+\beta|$,  
equation \eq{G0H} reduces to
\Eq{*}{
  G_0\Big(\frac{1}{2D}\ln\frac{v}{w}\Big)
  =C\ln\frac{v+w}{2}-\frac{(C\ln v)+(C\ln w)}{2} 
  =C\ln\frac{v+w}{2\sqrt{vw}}=C\ln\frac{\sqrt{\frac{v}{w}}+\sqrt{\frac{w}{v}}}{2}
}
Sunstituting $u:=\frac{1}{2D}\ln\frac{v}{w}$, it follows that $\sqrt{\frac{v}{w}}=e^{Du}$, hence 
the above equation yields that
\Eq{*}{
G_0(u)=C\ln(\cosh(Du)).
}
That is, we got the solutions listed in \emph{(iv)}.

%----------Hiperbolikus(a)
\emph{Case 4. Assume that $\gamma>0$ and $|a|=|b|$.} In this case, due to the condition $ad\neq 
bc$, both of the parameters $a$ and $b$ are different from zero. In order to get a simplier 
calculation, using the well-known identities $\sinh(t)=\frac12(e^{t}-e^{-t})$ and 
$\cosh(t)=\frac12(e^{t}+e^{-t})$ for $t\in\R$, we rewrite $\frac1{\ell'}$ and $H'$ to the form
\Eq{*}{
\frac{1}{\ell'(x)}
&=a\sinh(\sqrt{\gamma}x)+b\cosh(\sqrt{\gamma}x)
=\frac{a+b}{2}e^{\sqrt{\gamma}x}+\frac{b-a}{2}e^{-\sqrt{\gamma}x}=a\sgn(ab)e^{\sgn(ab)\sqrt{ \gamma}x},
\qquad(x\in I)
}
and
\Eq{*}{
H'(x)
=
\frac{c\sinh(\sqrt{\gamma}x)+d\sinh(\sqrt{\gamma}x)}{a\sinh(\sqrt{\gamma}x)+b\sinh(\sqrt{\gamma}x)}
&=\frac{d+c}{2a\sgn(ab)}e^{(1-\sgn(ab))\sqrt{\gamma}x}+\frac{d-c}{2a\sgn(ab)}e^{-(1+\sgn{ab})\sqrt{\gamma}x}\\
&=\frac{\sgn(ab)d-c}{2a}e^{-2\sgn(ab)\sqrt{\gamma}x}+\frac{c+\sgn(ab)d}{2a},
\qquad(x\in I).
}
Thus, solving \eq{not}, we obtain that there exist real constants $B$ and $E$ such that, for all 
$x\in I$, we have
\Eq{*}{
\ell(x)=
-\frac{1}{a\sqrt{\gamma}}e^{-\sgn(ab)\sqrt{\gamma}x}+E} 
and
\Eq{*}{
H(x)
=\frac{\sgn(ab)c-d}{4a\sqrt{\gamma}}e^{-2\sgn(ab)\sqrt{\gamma}x}+\frac{c+\sgn(ab)d}{2a}x+B.
}
Now, define $A, C, D,$ and $\alpha$ by
\Eq{*}{
A:=\frac{c+\sgn(ab)d}{2a},\quad
C:=\frac{d-\sgn(ab)c}{2a\sqrt{\gamma}},\quad
D:=-\frac{a\sqrt{\gamma}}{2},\quad\text{and}\quad
\alpha:=-\sgn(ab)\sqrt{\gamma}.
}
Obviously, $\alpha\neq 0$, furthermore, in view of the condition $ad\neq bc$, the constants $C$ and 
$D$ are also different from zero. Using these notations, we obtained that the functions $\ell$ 
and $H$ are of the form
\Eq{*}{
\ell(x)=\frac{1}{2D}e^{\alpha x}+E\qquad\text{and}\qquad
H(x)=-\frac{C}{2}e^{2\alpha x}+Ax+B
}
on the interval $I$. Therefore, with the substitutions $v:=\alpha x$ and $w:=\alpha y$, the 
equation \eq{G0H} reduces to
\Eq{*}{
G_0\Big(\frac{1}{2D}(e^v-e^w)\Big)
=-\frac{C}{2}e^{v+w}+\frac{\frac{C}{2}e^{2v}+\frac{C}{2}e^{2w}}{2}
=\frac{C}{4}(e^{2v}-2e^{v+w}+e^{2w})=\frac{C}{4}(e^{v}-e^{w})^2.
}
Let $u:=\frac{1}{2D}(e^v-e^w)$. Then $e^v-e^w=2Du$ and
\Eq{*}{
G_0(u)=\frac{C}{4}(2Du)^2=C(Du)^2.
}
That is, we obtain the solutions listed in \emph{(iii)}.

%----------Hiperbolikus(b)
\emph{Case 5. Assume that $\gamma>0$ and $|a|>|b|$.} Then it follows that $a^2-b^2>0$ and $a\neq0$. 
Denote $\alpha:=\frac{\sqrt{\gamma}}{2}\neq 0$ and define $\beta\in\R$ by the equation
\Eq{*}{
\sinh(2\beta)=\frac{\sgn(a)b}{\sqrt{a^2-b^2}}.
}
Then $\cosh(2\beta)=\frac{\sgn(a)a}{\sqrt{a^2-b^2}}$, therefore the identity
\Eq{*}{
0\neq \frac{1}{\ell'(x)}
=a\sinh(\sqrt{\gamma}x)+b\cosh(\sqrt{\gamma}x)=\sgn(a)\sqrt{a^2-b^2}\sinh(2\alpha x+2\beta),
\qquad(x\in I)
}
holds, and we have that $0\notin \alpha I+\beta$ (that is, the third condition in \eq{dom} holds). Solving the 
differential equations in \eq{not}, with the notations 
$A:=\frac{ac-bd}{a^2-b^2}$, $C:=\frac{ad-bc}{\sqrt{\gamma}(a^2-b^2)}\neq 0$ and 
$D:=\frac{\sqrt{\gamma(a^2-b^2)}}{2\sgn(a)}\neq 0$, we get that there exist constants $B$ and $E$ 
such that
\Eq{*}{
\ell(x)=\frac1{2D}\ln|\tanh(\alpha x+\beta)|+E,\qquad\text{and}\qquad
H(x)=C\ln|\sinh(2\alpha x+2\beta)|+Ax+B
}
hold for all $x\in I$. Consequently, introducing the notations $v:=\alpha x+\beta$ and 
$w:=\alpha y+\beta$, the equation \eq{G0H} reduces to
\Eq{*}{
G_0\Big(\frac{1}{2D}\ln\frac{\tanh v}{\tanh w}\Big)
&=C\ln|\sinh(v+w)|-\frac{C\ln|\sinh(2v)|+C\ln|\sinh(2w)|}{2}
=C\ln\frac{|\sinh(v+w)|}{\sqrt{\sinh(2v)\sinh(2w)}}\\
&=C\ln\frac{|\sinh v\cosh w+\sinh w\cosh v|}{\sqrt{4\sinh v\cosh v\sinh w\cosh w}}
=C\ln\frac{\sqrt{\frac{\tanh v}{\tanh w}}+\sqrt{\frac{\tanh w}{\tanh v}}}{2}.
}
(Note that, due to the condition $0\notin \alpha I+\beta$, in the above computation we have that $\tanh v$ and 
$\tanh w$ as well as $\sinh(2v)$ and $\sinh(2w)$ have the same sign.) Let now $u:=\frac{1}{2D}\ln\frac{\tanh v}{\tanh w}$. 
Then $\sqrt{\frac{\tanh v}{\tanh w}}=e^{Du}$ and
\Eq{*}{
G_0(u)=C\ln(\cosh(Du)).
}
That is, we obtain the solutions listed in \emph{(vi)}.

%----------Hiperbolikus(c)
\emph{Case 6. Assume finally that $\gamma>0$ and $|a|<|b|$.} Then it follows that $b^2-a^2>0$ and $b\neq0$. 
Let $\alpha:=\frac{\sqrt{\gamma}}2\neq 0$ and define the parameter $\beta\in\R$ by the equation
\Eq{*}{
\sinh(2\beta)=\frac{\sgn(b)a}{\sqrt{b^2-a^2}}.
}
Then we have that $\cosh(2\beta)=\frac{\sgn(b)b}{\sqrt{b^2-a^2}}$ and therefore
\Eq{*}{
\frac{1}{\ell'(x)}
=a\sinh(\sqrt{\gamma}x)+b\cosh(\sqrt{\gamma}x)
=\sgn(b)\sqrt{b^2-a^2}\cosh(2\alpha x+2\beta),\qquad(x\in I).
}
In view of the identity above, solving the differential equations in \eq{not}, with the notations 
$A:=\frac{bd-ac}{b^2-a^2}$, $C:=\frac{cb-ad}{\sqrt{\gamma}(b^2-a^2)}$ and 
$D:=\frac{\sqrt{\gamma(b^2-a^2)}}{2\sgn(b)}$, we get that there exist constants $B$ and $E$ such 
that
\Eq{*}{
\ell(x)=\frac1{D}\arctan(\tanh(\alpha x+\beta))+E\qquad\text{and}\qquad
H(x)=C\ln(\cosh(2\alpha x+2\beta))+Ax+B
}
hold for all $x\in I$. By their definitions, we can also see that the parameters $C$ and $D$ are 
different from zero. 

In order to determine $G_0$, firstly, we are going to shape the expression in the argument of 
$G_0$. Since, for all $t\in\R$, we have that $-1<\tanh t<1$, hence $-\frac\pi4<\arctan(\tanh t)<\frac\pi4$.
Therefore, for all $x,y\in I$, 
\Eq{tan}{
  -\frac\pi2<\arctan(\tanh(\alpha x+\beta))-\arctan(\tanh(\alpha y+\beta))<\frac\pi2.
}
On the other hand, by the addition theorem of the tangent function, we obtain that
\Eq{*}{
\tan\big(\arctan(\tanh(\alpha x+\beta))-\arctan(\tanh(\alpha y+\beta))\big)
=\frac{\tanh(\alpha x+\beta)-\tanh(\alpha y+\beta)}{1+\tanh(\alpha x+\beta)\tanh(\alpha y+\beta)}.
}
By the inequalities of \eq{tan}, with the substitutions $v:=e^{\alpha x+\beta}$ and $w:=e^{\alpha 
y+\beta}$, it follows that
\Eq{*}{
\arctan(\tanh(\alpha x+\beta))-\arctan(\tanh(\alpha y+\beta))
&=\arctan\frac{\tanh(\alpha x+\beta)-\tanh(\alpha y+\beta)}
{1+\tanh(\alpha x+\beta)\tanh(\alpha y+\beta)}\\
&=\arctan\frac{\frac{v^2-1}{v^2+1}-\frac{w^2-1}{w^2+1}}{1+\frac{v^2-1}{v^2+1}\frac{w^2-1}{w^2+1}}
=\arctan\frac{v^2-w^2}{1+v^2w^2}.
}
Hence the equation \eq{G0H} reduces to
\Eq{*}{
G_0\Big(\frac{1}{D}&\arctan\frac{v^2-w^2}{1+v^2w^2}\Big)\\
&=C\ln\big(\cosh(\alpha x+\alpha y+2\beta)\big)
-\frac{C\ln(\cosh(2\alpha x+2\beta))+C\ln(\cosh(2\alpha y+2\beta))}{2}\\
&=C\ln\frac{\cosh(\alpha x+\beta+\alpha y+\beta)}
{\sqrt{\cosh(2\alpha x+2\beta)\cosh(2\alpha y+2\beta)}}
=C\ln\frac{1+v^2w^2}{\sqrt{(v^4+1)(w^4+1)}}
=C\ln\frac{1}{\sqrt{1+\big(\frac{v^2-w^2}{1+v^2w^2}\big)^2}}.
}
Substituting $u:=\frac1D\arctan\frac{v^2-w^2}{1+v^2w^2}$, we get that
\Eq{*}{
G_0(u)
=C\ln\frac{1}{\sqrt{1+\tan^2(Du)}}=C\ln(\cos(Du)).
}
Therefore, in this case, we get the solutions listed in \emph{(vi)}.
\end{proof}

\def\MR#1{}
%\bibliography{publ,funcequ,refs2}

\begin{thebibliography}{10}

\bibitem{Acz47b}
J.~Aczél.
\newblock {The notion of mean values}.
\newblock {\em Norske Vid. Selsk. Forh., Trondhjem}, 19(23):83–86, 1947.

\bibitem{Acz48a}
J.~Aczél.
\newblock {On mean values}.
\newblock {\em Bull. Amer. Math. Soc.}, 54:392–400, 1948.

\bibitem{Acz48b}
J.~Aczél.
\newblock {On mean values and operations defined for two variables}.
\newblock {\em Norske Vid. Selsk. Forh., Trondhjem}, 20(10):37–40, 1948.

\bibitem{Acz56c}
J.~Aczél.
\newblock {On the theory of means}.
\newblock {\em Colloq. Math.}, 4:33–55, 1956.

\bibitem{AczKuc89}
J.~Aczél and M.~Kuczma.
\newblock {On two mean value properties and functional equations associated
  with them}.
\newblock {\em Aequationes Math.}, 38(2-3):216–235, 1989.

\bibitem{AczMakPal99}
J.~Aczél, Gy. Maksa, and Zs. Páles.
\newblock {Solution to a functional equation arising from different ways of
  measuring utility}.
\newblock {\em J. Math. Anal. Appl.}, 233(2):740–748, 1999.

\bibitem{AczMakPal01}
J.~Aczél, Gy. Maksa, and Zs. Páles.
\newblock {Solution of a functional equation arising in an axiomatization of
  the utility of binary gambles}.
\newblock {\em Proc. Amer. Math. Soc.}, 129(2):483–493, 2001.

\bibitem{Baj58}
M.~Bajraktarević.
\newblock {Sur une équation fonctionnelle aux valeurs moyennes}.
\newblock {\em Glasnik Mat.-Fiz. Astronom. Društvo Mat. Fiz. Hrvatske Ser.
  II}, 13:243–248, 1958.

\bibitem{BajPal09a}
Sz. Baják and Zs. Páles.
\newblock {Invariance equation for generalized quasi-arithmetic means}.
\newblock {\em Aequationes Math.}, 77:133–145, 2009.

\bibitem{Ber00}
L.~R. Berrone.
\newblock {A localization principle for classes of means}.
\newblock {\em Demonstratio Math.}, 33(3):557–566, 2000.

\bibitem{BerLom01}
L.~R. Berrone and A.~L. Lombardi.
\newblock {A note on equivalence of means}.
\newblock {\em Publ. Math. Debrecen}, 58(1-2):49–56, 2001.

\bibitem{DarMakPal04}
Z.~Daróczy, Gy. Maksa, and Zs. Páles.
\newblock {On two-variable means with variable weights}.
\newblock {\em Aequationes Math.}, 67(1-2):154–159, 2004.

\bibitem{DarPal01a}
Z.~Daróczy and Zs. Páles.
\newblock {On means that are both quasi-arithmetic and conjugate arithmetic}.
\newblock {\em Acta Math. Hungar.}, 90(4):271–282, 2001.

\bibitem{DarPal02c}
Z.~Daróczy and Zs. Páles.
\newblock {Gauss-composition of means and the solution of the
  {M}atkowski–{S}utô problem}.
\newblock {\em Publ. Math. Debrecen}, 61(1-2):157–218, 2002.

\bibitem{Def31}
B.~de~Finetti.
\newblock {{S}ul concetto di media}.
\newblock {\em Giornale dell' Instituto, Italiano degli Attuarii}, 2:369–396,
  1931.

\bibitem{GilPal01a}
A.~Gilányi and Zs. Páles.
\newblock {A regularity theorem for composite functional equations}.
\newblock {\em Arch. Math. (Basel)}, 77(4):317–322, 2001.

\bibitem{Jar07}
J.~Jarczyk.
\newblock {Invariance of weighted quasi-arithmetic means with continuous
  generators}.
\newblock {\em Publ. Math. Debrecen}, 71(3-4):279–294, 2007.

\bibitem{Jar09b}
J.~Jarczyk.
\newblock {Invariance of quasi-arithmetic means with function weights}.
\newblock {\em J. Math. Anal. Appl.}, 353(1):134–140, 2009.

\bibitem{Jar09a}
J.~Jarczyk.
\newblock {Regularity theorem for a functional equation involving means}.
\newblock {\em Publ. Math. Debrecen}, 75(1-2):123–135, 2009.

\bibitem{Jar10a}
J.~Jarczyk.
\newblock {Invariance in a class of {B}ajraktarević means}.
\newblock {\em Nonlinear Anal.}, 72(5):2608–2619, 2010.

\bibitem{Jar10b}
J.~Jarczyk.
\newblock On an equation involving weighted quasi-arithmetic means.
\newblock {\em Acta Math. Hungar.}, 129(1-2):96--111, 2010.

\bibitem{Jar13c}
J.~Jarczyk.
\newblock Determination of conjugate means by reducing to the generalized
  {M}atkowski-{S}ut\^o equation.
\newblock {\em Acta Math. Hungar.}, 139(1-2):1--10, 2013.

\bibitem{JarMat06}
J.~Jarczyk and J.~Matkowski.
\newblock {Invariance in the class of weighted quasi-arithmetic means}.
\newblock {\em Ann. Polon. Math.}, 88(1):39–51, 2006.

\bibitem{JarMat10}
W.~Jarczyk and J.~Matkowski.
\newblock {Embeddability of mean-type mappings in a continuous iteration
  semigroup}.
\newblock {\em Nonlinear Anal.}, 72(5):2580–2591, 2010.

\bibitem{JarMakPal04}
A.~Járai, Gy. Maksa, and Zs. Páles.
\newblock {On {C}auchy-differences that are also quasisums}.
\newblock {\em Publ. Math. Debrecen}, 65:381–398, 2004.

\bibitem{KahMat16}
P.~Kahlig and J.~Matkowski.
\newblock {Invariant means related to classical weighted means}.
\newblock {\em Publ. Math. Debrecen}, 89(3):373–387, 2016.

\bibitem{Kol30}
A.~N. Kolmogorov.
\newblock {{S}ur la notion de la moyenne}.
\newblock {\em Rend. Accad. dei Lincei (6)}, 12:388–391, 1930.

\bibitem{LeoMatTri16}
P.~Leonetti, J.~Matkowski, and S.~Tringali.
\newblock On the commutation of generalized means on probability spaces.
\newblock {\em Indag. Math. (N.S.)}, 27(4):945--953, 2016.

\bibitem{Los99}
L.~Losonczi.
\newblock {Equality of two variable weighted means: reduction to differential
  equations}.
\newblock {\em Aequationes Math.}, 58(3):223–241, 1999.

\bibitem{Los00a}
L.~Losonczi.
\newblock {Equality of {C}auchy mean values}.
\newblock {\em Publ. Math. Debrecen}, 57:217–230, 2000.

\bibitem{Los03a}
L.~Losonczi.
\newblock {Equality of two variable {C}auchy mean values}.
\newblock {\em Aequationes Math.}, 65(1-2):61–81, 2003.

\bibitem{Los06b}
L.~Losonczi.
\newblock {Equality of two variable means revisited}.
\newblock {\em Aequationes Math.}, 71(3):228–245, 2006.

\bibitem{LosPal11a}
L.~Losonczi and Zs. Páles.
\newblock {Equality of two-variable functional means generated by different
  measures}.
\newblock {\em Aequationes Math.}, 81(1-2):31–53, 2011.

\bibitem{MakPal08}
Z.~Makó and Zs. Páles.
\newblock {On the equality of generalized quasiarithmetic means}.
\newblock {\em Publ. Math. Debrecen}, 72:407–440, 2008.

\bibitem{Mat99a}
J.~Matkowski.
\newblock {Invariant and complementary quasi-arithmetic means}.
\newblock {\em Aequationes Math.}, 57(1):87–107, 1999.

\bibitem{Mat04c}
J.~Matkowski.
\newblock {Solution of a regularity problem in equality of {C}auchy means}.
\newblock {\em Publ. Math. Debrecen}, 64(3-4):391–400, 2004.

\bibitem{Mat11d}
J.~Matkowski.
\newblock {A functional equation related to an equality of means problem}.
\newblock {\em Colloq. Math.}, 122(2):289–298, 2011.

\bibitem{Mat14}
J.~Matkowski.
\newblock {Invariance identity in the class of generalized quasiarithmetic
  means}.
\newblock {\em Colloq. Math.}, 137(2):221–228, 2014.

\bibitem{MatNowWit17}
J.~Matkowski, M.~Nowicka, and A.~Witkowski.
\newblock {Explicit solutions of the invariance equation for means}.
\newblock {\em Results Math.}, 71(1-2):397–410, 2017.

\bibitem{MatPal15}
J.~Matkowski and Zs. Páles.
\newblock {Characterization of generalized quasi-arithmetic means}.
\newblock {\em Acta Sci. Math. (Szeged)}, 81(3-4):447–456, 2015.

\bibitem{MunMakMok99}
Á. Münnich, Gy. Maksa, and R.~J. Mokken.
\newblock {Collective judgement: combining individual value judgements}.
\newblock {\em Math. Social Sci.}, 37(3):211–233, 1999.

\bibitem{MunMakMok00}
Á. Münnich, Gy. Maksa, and R.~J. Mokken.
\newblock {$n$-variable bisection}.
\newblock {\em J. Math. Psych.}, 44(4):569–581, 2000.

\bibitem{Nag30}
M.~Nagumo.
\newblock {Über eine {K}lasse der {M}ittelwerte}.
\newblock {\em Japanese J. Math.}, 7:71–79, 1930.

\bibitem{Pal82a}
Zs. Páles.
\newblock {Characterization of quasideviation means}.
\newblock {\em Acta Math. Acad. Sci. Hungar.}, 40(3-4):243–260, 1982.

\bibitem{Pal84b}
Zs. Páles.
\newblock {On the characterization of means defined on a linear space}.
\newblock {\em Publ. Math. Debrecen}, 31(1-2):19–27, 1984.

\bibitem{Pal87d}
Zs. Páles.
\newblock {On the characterization of quasi-arithmetic means with weight
  function}.
\newblock {\em Aequationes Math.}, 32(2-3):171–194, 1987.

\bibitem{Pal02a}
Zs. Páles.
\newblock {Problems in the regularity theory of functional equations}.
\newblock {\em Aequationes Math.}, 63(1-2):1–17, 2002.

\bibitem{Pal03b}
Zs. Páles.
\newblock {A regularity theorem for composite functional equations}.
\newblock {\em Acta Sci. Math. (Szeged)}, 69:591–604, 2003.

\bibitem{Pal11}
Zs. Páles.
\newblock {On the equality of quasi-arithmetic and {L}agrangian means}.
\newblock {\em J. Math. Anal. Appl.}, 382(1):86–96, 2011.

\end{thebibliography}
%\bibliographystyle{plain}

\end{document}